\documentclass{article}%
\pdfoutput=1
\usepackage{graphicx}
\usepackage{amsmath, amsthm}
\usepackage{amsfonts}

\newtheorem{prop}{Proposition}[section]

\newtheorem{rem}{Remark}[section]

\newcommand{\ignore}[1]{}
\ignore{The text written between the parentheses here is ignored by latex}

\begin{document}

 \pagestyle{myheadings}
 \markboth{The complexity of the envelope of line and plane arrangements} {David Bremner, Antoine  Deza and Feng Xie}

\author{David Bremner, Antoine  Deza and Feng Xie}
\title{The complexity of the envelope\\ of line and plane arrangements}
\date{September 15, 2007}
  \maketitle

\begin{abstract}
A facet of an hyperplane arrangement is called external if it belongs to exactly one 
bounded cell. The set of all external facets forms the envelope of the arrangement.
The number of external facets of a simple 
arrangement defined by $n$ hyperplanes in dimension $d$ is hypothesized to be at 
least $d{n-2 \choose d-1}$. In this note we show that, for simple arrangements of $4$ lines or more, 
the minimum number of external facets is equal to $2(n-1)$, and for simple arrangements of $5$ planes
or more, the minimum number of external facets is between $\frac{n(n-2)+6}{3}$ and $(n-4)(2n-3)+5$.
\end{abstract}

\section{Introduction}
Let $\mathcal{A}_{\:d,n}$ be a simple arrangement formed by
$n$ hyperplanes in dimension $d$. We recall that an arrangement is
called simple if $n\geq d+1$ and any  $d$ hyperplanes intersect at a
distinct point. The closures of connected components of the complement of the hyperplanes 
forming ${\mathcal A}_{\:d,n}$ are called the cells, or $d$-faces, of the arrangement.
For $k=0,\dots,d-1$, the $k$-faces of ${\mathcal A}_{\:d,n}$ are the $k$-faces of its cells. 
A facet is a $(d-1)$-face of $\mathcal{A}_{\:d,n}$, and a  facet  belonging to exactly one bounded
cell is called an external facet. Equivalently, an external facet is a bounded facet which
belongs to an unbounded cell. For $k=0,\dots,d-2$, an external $k$-face is a $k$-face
belonging to an external facet. Let $f_k^0({\mathcal A}_{\:d,n})$ denote the number of 
external $k$-faces of ${\mathcal A}_{\:d,n}$. The set of all external facets forms the 
envelope of the arrangement. It was hypothesized
in~\cite{dtz06} that any simple arrangement $\mathcal{A}_{\:d,n}$ has at least 
$d{n-2 \choose d-1}$ external facets. In Section~\ref{2d}, we show that a 
simple arrangement of $n$ lines has at least $2(n-1)$ external facets for $n\geq 4$,
and that this bound is tight. In section~\ref{3d}, we show that  a 
simple arrangement of $n$ planes has at least $\frac{n(n-2)+6}{3}$ external facets
for $n\geq 5$, and exhibit a simple plane arrangement with $(n-4)(2n-3)+5$ external facets.
For polytopes and arrangements, we refer to the books of Edelsbrunner~\cite{E87},
Gr\"unbaum~\cite{G03} and Ziegler~\cite{Z95} and the references therein.

\section{The complexity of the envelope of line arrangements} \label{2d}
\subsection{A lower bound}
\begin{figure}[htb]
\begin{center}
\includegraphics[height=2.5cm]{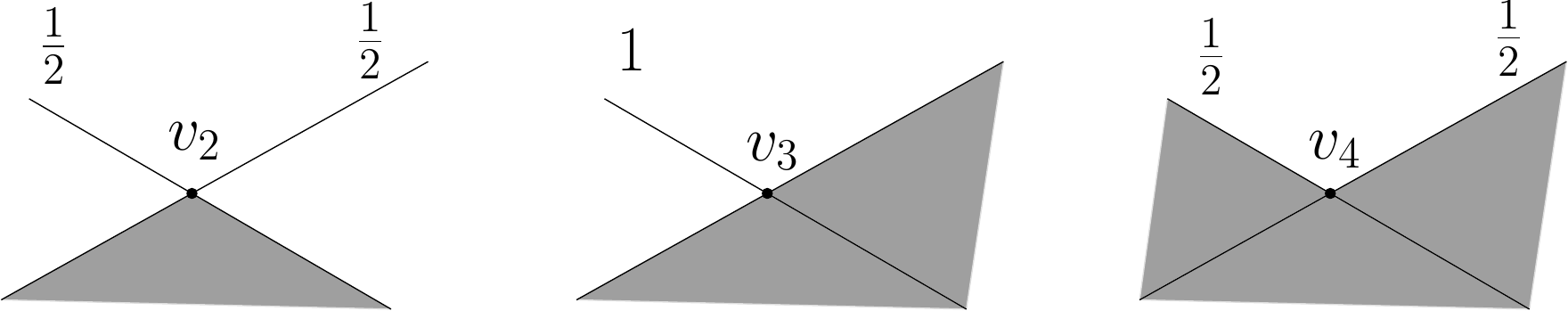}
\caption{The weight distribution for the lines of an arrangement (the shaded area 
corresponds to the bounded cells).} \label{Fig_Weight_Distribution}
\end{center}
\end{figure}

\begin{prop}\label{ext2dLB}
For $n\geq 4$, 
a simple line arrangement has at least $2(n-1)$ external facets.
\end{prop}
\begin{proof}
The external vertices of a line arrangement can be
divided into three types, namely $v_2$, $v_3$ and $v_4$, corresponding to external vertices
respectively incident to 2, 3, and 4 bounded edges. 
Let us assign to each external vertex $v$ a weight of 1 and redistribute it to 
the 2 lines intersecting at $v$ the following way: If $v$ is incident to
exactly 1 unbounded edge, then give weight 1 to the line containing 
this edge, and weight 0 to the other line containing $v$; if $v$ is incident
to $2$ or $0$ unbounded edges,  then give weight $0.5$ to each of the
$2$ lines intersecting at $v$. See Figure~\ref{Fig_Weight_Distribution}
for an illustration of the weight distribution.
A total of $f^0_0({\mathcal A}_{\:2,n})$ weights
is distributed and we can also count this quantity line-wise.
The end vertices of a line being of type $v_2$ or $v_3$, we have
three types of lines, $h_{2,2}, h_{2,3}$ and $h_{3,3}$, 
according to the possible types of their end-vertices. 
 As a line of type $h_{3,3}$ contains $2$ vertices of type $v_3$, its weight
is at least $2$.  Similarly the weight of a line of type $h_{2,2}$  weight
is at least $1$. Remarking that a line of type $h_{2,3}$ contains at 
least one vertex of type $v_4$ yields that the weight of a line of 
type $h_{2,3}$ is at least $1+0.5+0.5=2$. For $n\geq 4$ the number 
of lines of type $h_{2,2}$ is at most 2 as otherwise the envelope
would be convex which is impossible, see for example~\cite{EGT96}.
Therefore, counting the total distributed weight line-wise, we have 
$f^0_0({\mathcal A}_{\:2,n})\geq 2n-2$. Since for a line arrangement
the number of external facets $f^0_1({\mathcal A}_{\:2,n})$
is equal to the number of external vertices $f^0_0({\mathcal
A}_{\:2,n})$, we have $f^0_1({\mathcal A}_{\:2,n})\geq 2(n-1)$.
\end{proof}

\subsection{A line arrangement attaining the lower bound}
\begin{figure}[hbt]
\begin{center}
\includegraphics[height=10cm]{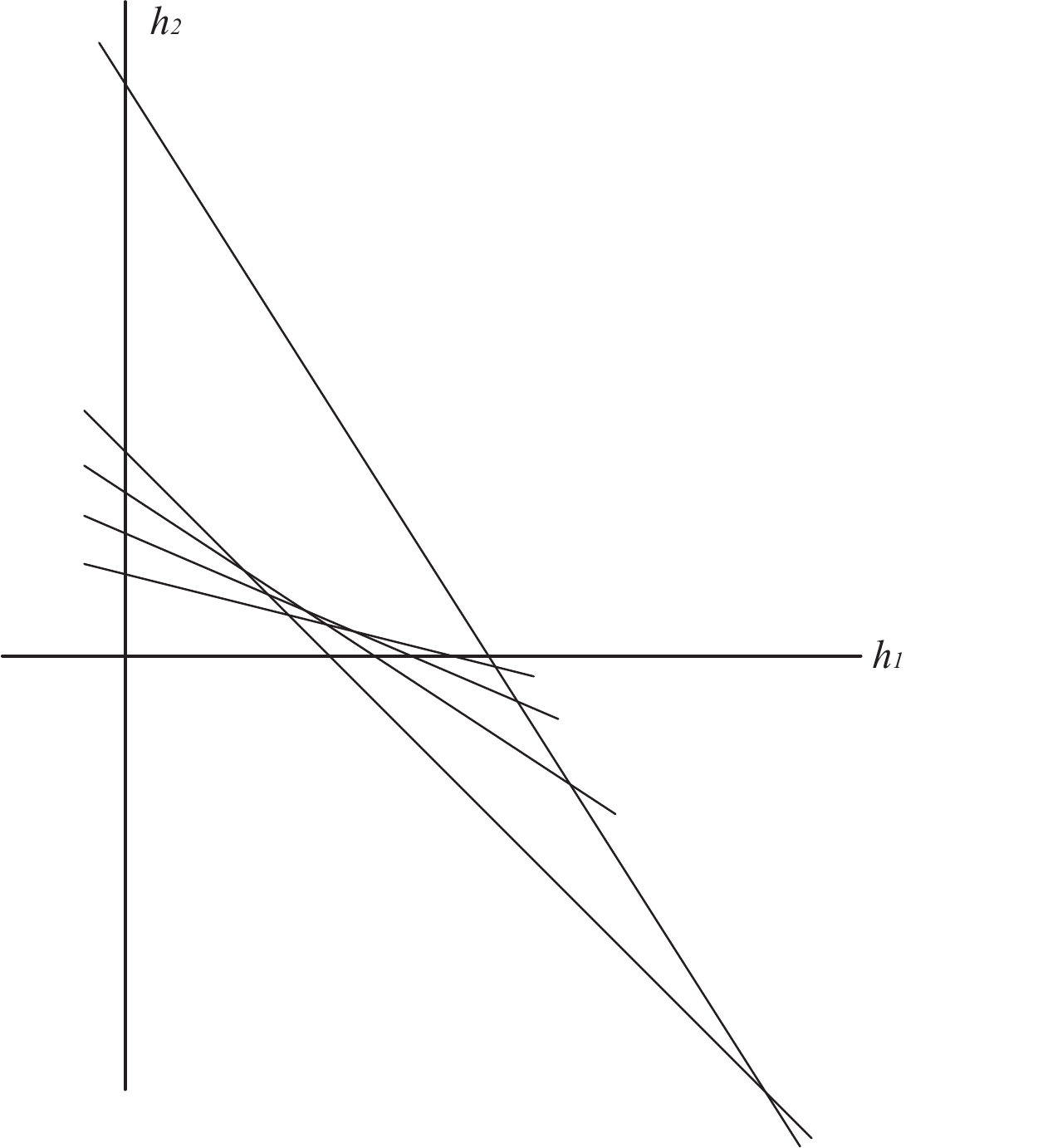}
 \caption{An arrangement
combinatorially equivalent to $\mathcal{A}^o_{\:2,7}$} \label{A072}
\end{center}
\end{figure}

For $n\geq 4$,  consider the following simple line arrangement:
$\mathcal{A}^o_{\:2,n}$ is
 made of the $2$ lines $h_1$ and $h_2$
forming, respectively, the $x_1$ and $x_2$ axis, and $(n-2)$ lines
defined by their intersections with $h_1$ and $h_2$. We have
$h_k\cap h_1=\{1+(k-3)\varepsilon,0\}$ and $h_k\cap
h_2=\{0,1-(k-3)\varepsilon\}$ for $k=3,4,\dots,n-1$, and $h_n\cap
h_1=\{2,0\}$ and $h_n\cap h_1=\{0,2+\varepsilon\}$ where
$\varepsilon$ is a constant satisfying
$0<\varepsilon<1/(n-3)$. See Figure~\ref{A072} for an
arrangement combinatorially equivalent to  $\mathcal{A}^o_{\:2,7}$.
One can easily check that  $\mathcal{A}^o_{\:2,7}$ has $2(n-1)$ external facets
and therefore the lower bound given in Proposition~\ref{ext2dLB}
is tight.

\begin{prop}
For $n\geq 4$, the minimum possible number of external facets of a simple  line
arrangement is $2(n-1)$.
\end{prop}

\section{The complexity of the envelope of plane arrangements} \label{3d}
\subsection{A lower bound}
\begin{prop}\label{ext3dLB}
For $n\geq 5$, 
a simple plane arrangement has at least $\frac{n(n-2)+6}{3}$ external facets.
\end{prop}
\begin{proof}
Let $h_i$ for $i=1,2,\dots,n$ be the planes forming the arrangement $\mathcal{A}_{\:3,n}$.
For $i=1,2,\dots,n$, the external vertices of the line arrangement $\mathcal{A}_{\:3,n} \cap h_i$ 
are external vertices of the plane arrangement $\mathcal{A}_{\:3,n}$.
For $n\geq 5$, the line arrangement $\mathcal{A}_{\:3,n} \cap h_i$ has at least $2(n-2)$ external facets
by Proposition~\ref{ext2dLB}, i.e., at least $2(n-2)$ external
vertices. Since an external vertex of  $\mathcal{A}_{\:3,n}$ belongs to 3 planes, it is counted three
times. In other words, the number of external vertices of  $\mathcal{A}_{\:3,n}$ satisfies
$f^0_0({\mathcal A}_{\:3,n}) \geq \frac{2n(n-2)}{3}$ for $n\geq 5$. As the 
union of all of the bounded cells is a piecewise linear ball, see~\cite{D06}, 
the Euler characteristic of the boundary gives 
$f^0_0({\mathcal A}_{\:3,n}) - f^0_1({\mathcal A}_{\:3,n}) +
f^0_2({\mathcal A}_{\:3,n}) = 2$.  Since an external vertex belong to at least 
3 external edges, we have $2 f^0_1({\mathcal A}_{\:3,n}) \geq 3 f^0_0({\mathcal A}_{\:3,n})$.
Thus, we have
$2f^0_2({\mathcal A}_{\:3,n}) \geq f^0_0({\mathcal A}_{\:3,n}) + 4$. As $f^0_0({\mathcal A}_{\:3,n}) \geq \frac{2n(n-2)}{3}$,
it gives $f^0_2({\mathcal A}_{\:3,n}) \geq \frac{n(n-2)+6}{3}$
\end{proof}

\subsection{A plane arrangement with few external facets}
\begin{figure}[htb]
\begin{center}
\includegraphics[height=13cm]{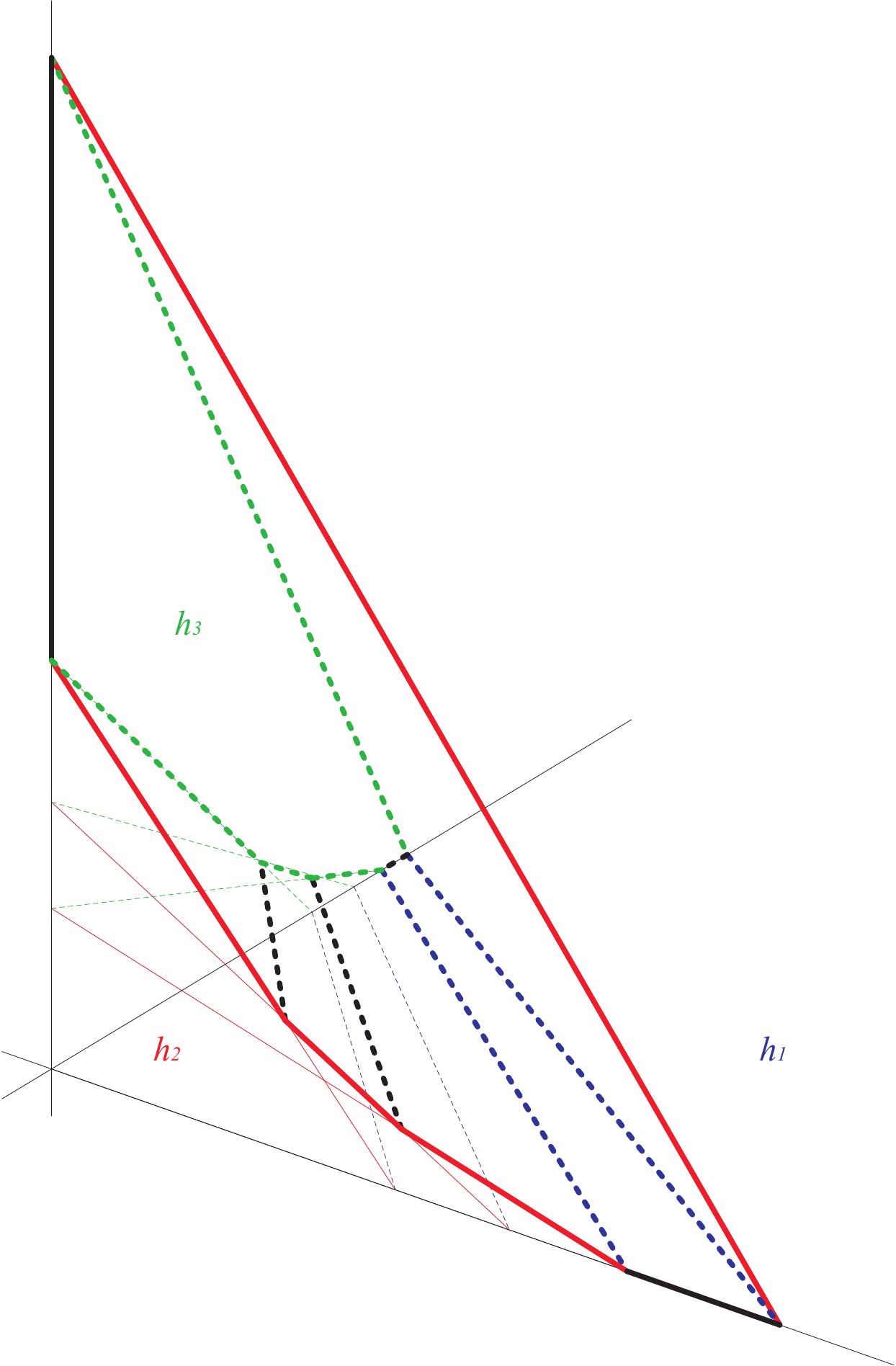} \caption{An arrangement
combinatorially equivalent to $\mathcal{A}^o_{\:3,7}$} \label{A073}
\end{center}
\end{figure}

For $n\geq 5$, we consider following simple plane arrangement:
$\mathcal{A}^o_{\:3,n}$ is made of the
$3$ planes $h_1$, $h_2$ and $h_3$ corresponding, respectively, to
$x_3=0$, $x_2=0$ and $x_1=0$, and $(n-3)$ planes defined by their
intersections with the $x_1$, $x_2$ and $x_3$ axis. We have $h_k\cap
h_1\cap h_2=\{1+2(k-4)\varepsilon,0,0\}$, $h_k\cap h_1\cap
h_3=\{0,1+(k-4)\varepsilon,0\}$ and $h_k\cap h_2\cap
h_3=\{0,0,1-(k-4)\varepsilon\}$ for $k=4,5,\dots,n-1$, and $h_n\cap
h_1\cap h_2=\{3,0,0\}$, $h_n\cap h_1\cap h_3=\{0,2,0\}$ and $h_n\cap
h_2\cap h_3=\{0,0,3+\varepsilon\}$ where $\varepsilon$ is a constant
satisfying $0<\varepsilon<1/(n-4)$. See Figure~\ref{A073} for
an illustration of an arrangement combinatorially equivalent to
$\mathcal{A}^o_{\:3,7}$ where, for clarity, only the bounded cells
belonging to the positive orthant are drawn.

\begin{figure}[htb]
\begin{center}
\includegraphics[height=11.3cm]{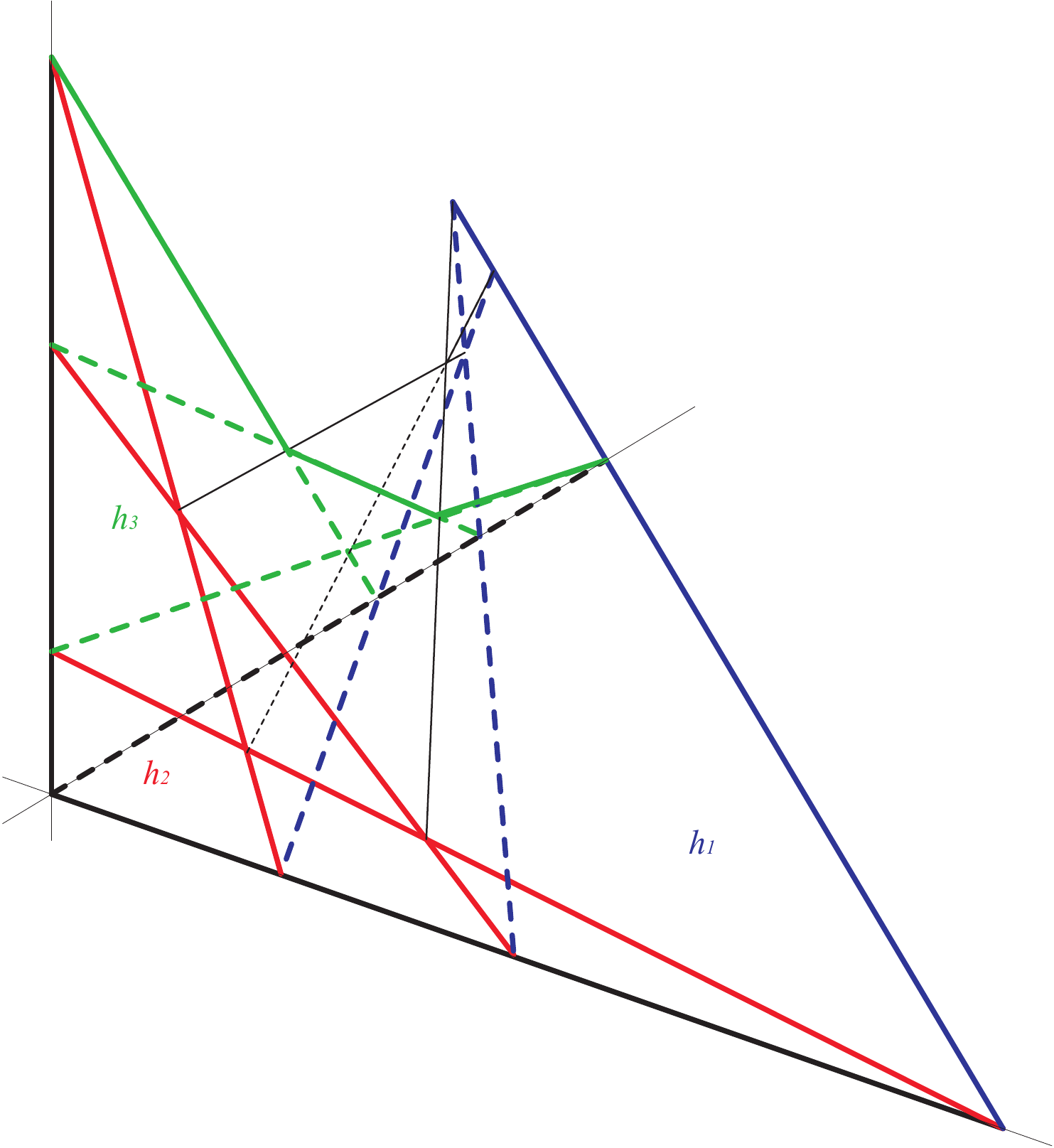} \caption{An arrangement
combinatorially equivalent to $\mathcal{A}^*_{\:3,6}$} \label{A63}
\end{center}
\end{figure}

We first check by induction that the arrangement $\mathcal{A}^*_{\:3,n}$ 
formed by the first $n$ planes of $\mathcal{A}^o_{\:n+1,3}$ has $2(n-2)(n-3)$ external facets. 
The arrangement $\mathcal{A}^*_{\:3,n}$  is combinatorially equivalent to the plane cyclic
arrangement which is dual to the cyclic polytope, see~\cite{FR01} for combinatorial properties of the 
(projective) cyclic arrangement in general dimension. 
See Figure~\ref{A63} for an illustration of $\mathcal{A}^*_{\:3,6}$. 
Let $H^+_3$ denote the half-space defined by $h_3$ and containing the positive orthant,
and  $H^-_3$ the other half-space defined by $h_3$. 
The union of the bounded cells of ${\mathcal A}^*_{3,n}$ in $H^-_3$ is 
combinatorially  equivalent to the bounded cells of ${\mathcal A}^*_{3,n-1}$ and therefore
has $2(n-3)(n-4)$ facets on its boundary by induction hypothesis, 
including $n-3\choose 2$ bounded facets contained in $h_3$.
These $n-3\choose 2$ bounded facets also belong to a bounded cell of
${\mathcal A}^*_{3,n}$ in $H^-_3$ and therefore are not external facets
of ${\mathcal A}^*_{3,n}$. Thus, 
the number of external facets  of ${\mathcal A}^*_{3,n}$ belonging to a bounded cell in $H^-_3$ is
$2(n-3)(n-4)-{n-3\choose 2}$.
The union of the bounded cells of ${\mathcal A}^*_{3,n}$ in $H^+_3$ can be viewed 
as a simplex cut by $n-4$ sliding down planes. It has $2{n-2 \choose 2}+2(n-3)=n(n-3)$ facets on its boundary, including the $n-3\choose 2$ bounded facets contained in $h_3$ belonging to a bounded cell of ${\mathcal A}^*_{3,n}$ 
in $H^-_3$. Thus, the number of external facets  of ${\mathcal A}^*_{3,n}$ belonging to a 
bounded cell in $H^+_3$ is $n(n-3)-{n-3\choose 2}$.
Therefore, ${\mathcal A}^*_{3,n}$ has $n(n-3)+2(n-3)(n-4)-2{n-3\choose 2}=2(n-2)(n-3)$
external facets. 
We now consider how the addition of $h_n$ to ${\mathcal A}^*_{3,n-1}$ impacts  the number of external facets. This impact is similar in nature to the addition of $h_n$ to the first $n-1$ lines of $\mathcal{A}^o_{\:2,n}$. The addition of $h_n$ creates $n \choose 2$ new bounded cells: 
one above $h_1$ that we call the \emph{$n$-shell}, and the other ones being below $h_1$. 
The $n$-shell turns $n-4$ external facets of ${\mathcal A}^*_{3,n-1}$ above $h_1$ into internal facets of
$\mathcal{A}^o_{\:3,n}$, and adds 3 external facets. For each external facet of ${\mathcal A}^*_{3,n-1}$ belonging to $h_1$ 
which is turned into an internal facet of $\mathcal{A}^o_{\:3,n}$,  one external facet of $\mathcal{A}^o_{\:3,n}$ on $h_n$  
and not incident to $h_1$ is added. Below $h_1$, the addition of $h_n$ creates $3(n-4)+2$ new external facets of $\mathcal{A}^o_{\:3,n}$
with an edge on $h_1$. Finally, $n-4$ new external facets belonging to $h_1$ and bounded by $h_n$ are created from 
unbounded facets of ${\mathcal A}^*_{3,n-1}$. Thus, the total number of external facets of ${\mathcal A}^o_{3,n}$ 
is $2(n-3)(n-4)-(n-4)+3+(3(n-4)+2)+(n-4)=(n-4)(2n-3)+5$.

\begin{figure}[htb]
\begin{center}
\includegraphics[height=6.3cm]{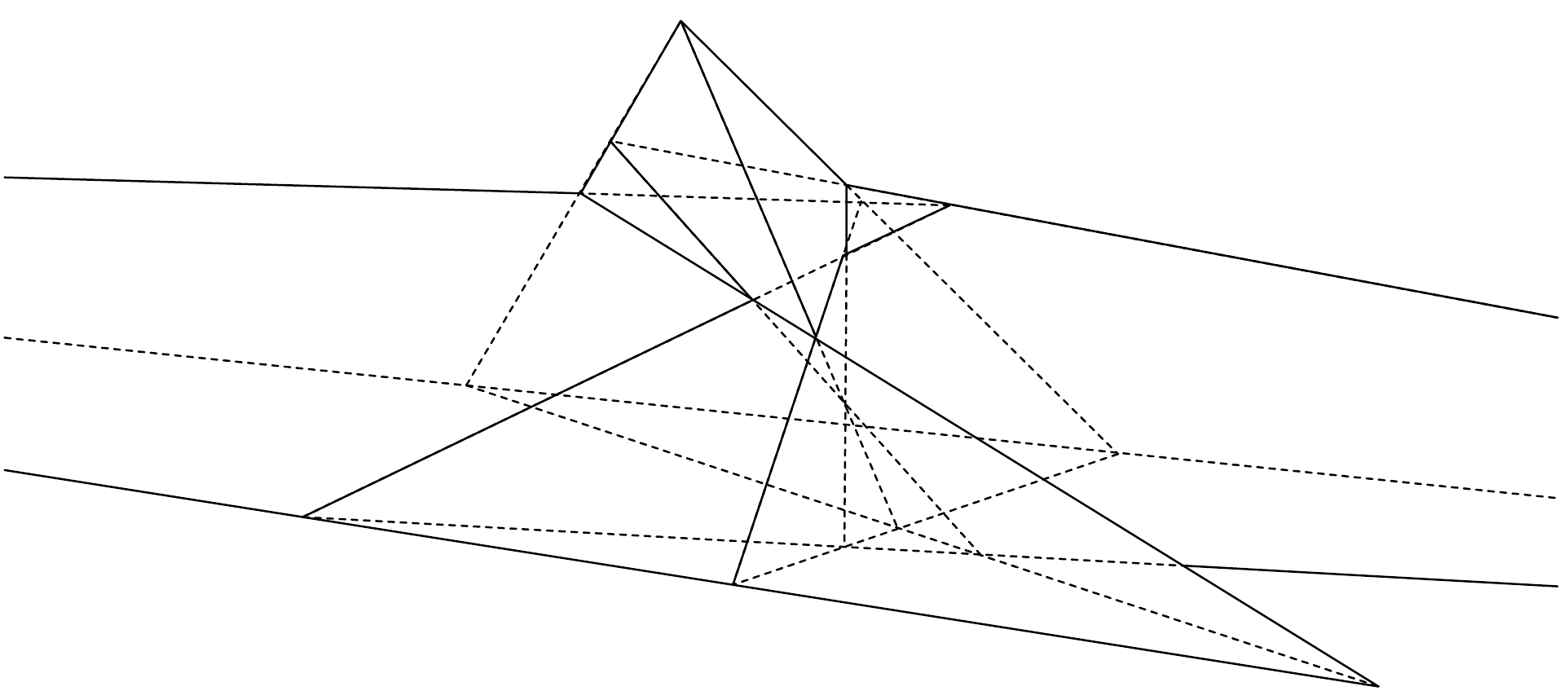} \caption{An arrangement
formed by 6 planes and having $22$ external facets} \label{A63_29}
\end{center}
\end{figure}

\begin{rem}
We do not believe that  ${\mathcal A}^{o}_{\:3,n}$ minimizes the 
number of external facets. Among the 43 simple combinatorial types of
arrangements formed by 6 planes, the minimum number of external facets is
$22$ while $\mathcal{A}^o_{\:3,6}$  has 23 external facets. See Figure~\ref{A63_29}
for an illustration of the combinatorial type of one of the two
simple arrangements with $6$ planes having 22  external facets.
The far away vertex on the right and 3 bounded edges incident to it 
are cut off  (same for the far away vertex on the left) so the 10 bounded cells 
of the arrangement appear not too small.
 \end{rem}

$\:$\\
\noindent{\bf Acknowledgments $\:$} Research supported by NSERC Discovery grants, by 
MITACS grants, by the Canada Research Chair program, and by the
Alexander von Humboldt Foundation.

\newpage
\noindent {\small David Bremner}\\
{\sc Faculty of Computer Science,}\\
{\sc  University of New Brunswick, New Brunswick, Canada}. \\
{\em Email}: bremner{\small @}unb.ca\\\\

\noindent {\small Antoine Deza, Feng Xie}\\
{\sc Department of Computing and Software,}\\
{\sc  McMaster University, Hamilton, Ontario, Canada}. \\
{\em Email}: deza, xief{\small @}mcmaster.ca

\end{document}